\centerline{MATRIX-VARIATE GROWTH-DECAY MODELS}
\vskip.5cm
\noindent
\centerline{A.M. Mathai}
\centerline {Department of Mathematics and Statistics}
\centerline{McGill University}
\centerline{805 Sherbrooke Street West, Montreal, Canada, H3A 2K6}
\vskip1cm
\noindent
\hrule
\vskip.5cm
\noindent
ABSTRACT
\vskip.5cm
Input-output, growth-decay, production-consumption type situations abound
 in many practical problems. When the input and output variables are
 independently gamma distributed, various aspects of the residual effect are
 already tackled by the author. Matrix-variate analogues, their connections
 to quadratic and bilinear forms, matrix-variate Whittaker functions, and many
 properties of such matrix-variate functions, in the real as well as complex
 case, are discussed here.
\vskip.3cm
\noindent
\hrule
\vskip.5cm
\noindent
1.{\hskip.5cm}INTRODUCTION
\vskip.5cm
If the input variable $x_1$ and the output variable $x_2$ are independently
 gamma distributed then it is known that the residual effect $y=x_1-x_2$ has a
 density which can be expressed in terms of Whittaker functions in the scalar
 case. A particular case is a Laplace density.  When $x_j$ has the density
$$f_j(x_j)={{x_j^{\alpha_j-1}{\rm
e}^{-x_j/{\beta_j}}}\over{\beta_j^{\alpha_j}\Gamma(\alpha_j)}},~x_j>0,~{\rm
Re}(\alpha_j)>0,~{\rm Re}(\beta_j)>0,~j=1,2\eqno(1.1)
$$and $f_j(x_j)=0$ elsewhere, where ${\rm Re}(\cdot)$ denotes the real part of
 $(\cdot)$, the density of $y=x_1-x_2$, denoted by $f(y)$, is given by
$$f(y)=\cases{c_1~y^{{{\alpha_1+\alpha_2}\over2}-1}{\rm
e}^{-{{y}\over2}\left({{1}\over{\beta_1}}-{{1}\over{\beta_2}}\right)}\cr
\cr
\times
~W_{{{\alpha_1-\alpha_2}\over2},{{1-\alpha_1-\alpha_2}\over2}}(\beta_0 y),
~y>0\cr
\cr
c_2~(-y)^{{{\alpha_1+\alpha_2}\over2}-1}{\rm
e}^{{{y}\over2}\left({{1}\over{\beta_2}}-{{1}\over{\beta_1}}\right)}\cr
\cr
\times~W_{{{\alpha_2-\alpha_1}\over2},{{1-\alpha_1-\alpha_2}\over2}}(-\beta_
0y),~y<0\cr}\eqno(1.2)
$$where
$$\eqalignno{\beta_0&={{1}\over{\beta_1}}+{{1}\over{\beta_2}},\cr
c_1^{-1}&=\Gamma(\alpha_1)\beta_1^{(\alpha_1-\alpha_2)/2}\beta_2^{(\alpha_2-
\alpha_1)/2}(\beta_1+\beta_2)^{(\alpha_1+\alpha_2)/2},\cr
c_2^{-1}&=\Gamma(\alpha_2)\beta_1^{(\alpha_1-\alpha_2)/2}\beta_2^{(\alpha_2-
\alpha_1)/2}(\beta_1+\beta_2)^{(\alpha_1+\alpha_2)/2}\cr}
$$and  $W_{.,.}(\cdot)$ is a Whittaker function, see also Mathai (1993a).
\vskip.2cm
Independent sums of such residual effects are shown to be connected to the
 distributions of quadratic and bilinear forms in Gaussian random variables,
 see Mathai and Provost (1992),  Mathai, Provost and Hayakawa (1995) and
Rao (1973).
 Laplacianness of quadratic and bilinear forms or the conditions under which a
 quadratic or bilinear form is distributed as a Laplace variable or as a gamma
 difference is also investigated. Let $S$ be a $p\times p$ non-singular central
 Wishart matrix and let it be partitioned as
$$S=\left[\matrix{S_{11}&S_{12}\cr
S_{21}&S_{22}\cr}\right]\eqno(1.3)
$$where $S_{11}$ is $r\times r,~r<p$ then it is shown that the distribution of
 $S_{12}$ is associated with that of certain bilinear forms. Note that $S_{12}$
 can also be looked upon as a certain covariance structure when $S$ is a
 corrected sample sum of products matrix. A large number of results in this
 category may be seen from Mathai (1993b) and Mathai, Provost and Hayakawa
 (1995).
\vskip.2cm
Solar neutrinos are captured on earth and the data for over twenty years are
 available. The summarized data, by taking a five point moving average to smooth
 out local disturbances, show a peculiar pattern of a 9-year cyclic behavior
 where within each cycle the pattern is a slow rising and rapidly decreasing
 curve with one real peak and several small humps in between. This type of
 pattern is also seen in many other problems such as the production of
 melatonin in human body. In Haubold and Mathai (1995) it is shown that this
 type of curve could be generated by combining three residual variables where
 the input and output variables are of gamma type, and some heuristic
 interpretations of these variables are also given there.
\vskip.2cm
One matrix-variate analogue of gamma type input and gamma type output gives
 the residual variable $Y=X_1-X_2$ where $X_1$ and $X_2$ are independently
 distributed $p\times p$ matrix-variate real gamma random variables with the
 densities
$$g_j(X_j)={{|B_j|^{\alpha_j}}\over{\Gamma_p(\alpha_j)}}|X_j|^{\alpha_j-{{p+
1}\over2}}{\rm e}^{-{\rm tr}(B_jX_J)},\eqno(1.4)
$$for $X_j=X_j'>0,~B_j=B_j'>0,~{\rm Re}(\alpha_j)>{{p-1}\over2},$
 and $g_j(X_j)=0,~j=1,2$ elsewhere, where for example,
$$\Gamma_p(\alpha)=\pi^{p(p-1)/4}\Gamma(\alpha)\Gamma\left(\alpha-{1\over2}
\right)...\Gamma\left(\alpha -{{p-1}\over2}\right),~{\rm Re}(\alpha)>{{p-1}\over2}\eqno(1.5)
$$is the matrix-variate gamma function in the real case, with the integral
representation
$$\Gamma_p(\alpha)=\int_{X=X'>0}|X|^{\alpha -{{p+1}\over2}}{\rm e}^
{-{\rm tr}(X)}{\rm d}X.\eqno(1.6)
$$When $B_j={1\over2}A_j$ for some $A_j=A_j'>0$ one has the non-singular
 central Wishart density from (1.4). Then if $X_1$ and $X_2$ are two
 independent central Wishart matrices then the residual variable is
$$Y=X_1-X_2.
$$One can look at element-wise differences in $X_1-X_2$ and consider the
 configuration of residual variables or one can consider the difference in
 terms of the definiteness of the $p\times p$ matrices $X_1$ and $X_2$. Since
 $X_1$ and $X_2$ are real symmetric positive definite matrices in (1.4), $Y$
 could be positive definite or negative definite or indefinite or semidefinite.
 A special
 situation will be when $X_1$ and $X_2$ are oriented matrices such that $Y$ is
 either nonnegative definite or $-Y$ is nonnegative definite. In this case one
 can proceed in the following way to compute the density of $Y$.
The joint density of $X_1$ and $X_2$, denoted by
$f(X_1,X_2)$, is the product of $f_1(X_1)$ and
$f_2(X_2)$ due to statistical independence. That is,
$$
\eqalignno{f(X_1,X_2){\rm d}X_1{\rm
d}X_2&=\delta |X_1|^{\alpha_1-{{p+1}\over2}}|X_2|^{\alpha_2-{{p+1}\over2}}{\rm
e}^{-{\rm tr}(B_1X_1+B_2X_2)}{\rm d}X_1{\rm d}X_2&(1.7)\cr
\noalign{\hbox{where}}
\delta
&={{|B_1|^{\alpha_1}|B_2|^{\alpha_2}}\over{\Gamma_{p}(\alpha_1)\Gamma_{p}
(\alpha_2)}}.&(1.8)\cr}
$$Consider the transformation $U=X_1-X_2$ and $V=X_1$.
The Jacobian is unity and the joint density of $U$ and
$V$, denoted by $g(U,V)$, is given by
$$g(U,V){\rm d}U{\rm
d}V=\delta |V|^{\alpha_1-{{p+1}\over2}}|V-U|^{\alpha_2-{{p+1}\over2}}{\rm
e}^{-{\rm tr}(B_1V+B_2(V-U))}{\rm d}U{\rm d}V.\eqno(1.9)
$$The density of $U$, denoted by $g(U)$, is available by
integrating out $V$ from (1.9). Note that when $U$ is
oriented as described earlier there are only two possible
situations, $U>0,V>U$ and $U<0,V>0$. Let us consider
these two situations separately. For $U>0, V>U$ let
$$
\eqalignno{g_1(U)&=\int_{V>U}g(U,V){\rm d}V\cr
&=\delta\int_{V>U}|V|^{\alpha_1-{{p+1}\over2}}|V-U|^{\alpha_2-{{p+1}\over2}}{\rm
e}^{-{\rm tr}(B_1V+B_2(V-U))}{\rm d}V\cr
&=\delta|U|^{\alpha_1+\alpha_2 -{{p+1}\over2}}{\rm e}^{-{\rm
tr}(B_1U)}\int_{Z>0}|Z|^{\alpha_2-{{p+1}\over2}}\cr
&\times\>|I+Z|^{\alpha_1-{{p+1}\over2}}{\rm e}^{-{\rm
tr}\left(U^{1\over2}(B_1+B_2)U^{1\over2}Z\right)}{\rm
d}Z&(1.10)\cr}
$$by making the transformations $Y=V-U$ and
$Z=U^{-{1\over2}}YU^{-{1\over2}}$ for fixed $U$. The integral in (1.10) even
 when
 $p=1$, that is in the scalar case, is difficult. In the scalar case it can be
 evaluated in terms of a Whittaker function which can then be represented as
 a linear combination of two confluent hypergeometric functions at least in some
 special cases. But when $p>1$ no such representation is available due to the
 fact that
when $p>1$
$$\int_{X=X'>0}\ne \int_{0<X=X'<1}+\int_{X=X'\ge I}
$$because $0<X<I$ implies only that the eigenvalues of $X$ are between zero and
 one and $X\ge I$ gives the eigenvalues greater than or equal to one. But a
 whole class of matrices with some eigenvalues  less than one and others
 greater than or equal to one are left out from the set $X=X'>0$. Hence the
 integral in (1.10) is not available even from a matrix-variate version of a
 confluent hypergeometric function. We need the Whittaker function of matrix
 argument to deal with this integral. This will be defined next.
\vskip.5cm
\noindent
2.{\hskip.5cm}WHITTAKER FUNCTION OF MATRIX ARGUMENT
\vskip.5cm

When studying some properties of the matrix-variate Laplace transform, as
 applied to special functions of matrix argument, it was necessary to come up
 with an expression for an integral of the following type.
$$
\eqalignno{\int_{X>A}|X|^{\alpha -{{p+1}\over2}}{\rm
e}^{-{\rm tr}(X)}{\rm d}X&={\rm e}^{-{\rm
tr}(A)}\int_{Y>0}|Y+A|^{\alpha -{{p+1}\over2}}{\rm
e}^{-{\rm tr}(Y)}{\rm d}Y\cr
&=|A|^{\alpha}{\rm e}^{-{\rm
tr}(A)}\int_{Z>0}|I+Z|^{\alpha -{{p+1}\over2}}{\rm
e}^{-{\rm tr}(AZ)}{\rm d}Z\cr}
$$by making the transformations $Y=X-A,
Z=A^{-{1\over2}}YA^{-{1\over2}}$. Hence Mathai and Pederzoli (1996) defined a
function $M(\alpha,\beta ;A)$ as follows:
$$M(\alpha,\beta ;A)=\int_{X=X'>0}|X|^{\alpha
-{{p+1}\over2}}|I+X|^{\beta -{{p+1}\over2}}{\rm
e}^{-{\rm tr}(AX)}{\rm d}X\eqno(2.1)
$$ for $A=A'>0, {\rm Re}(\alpha )>{{p-1}\over2}$ and no
restrictions on $\beta$. Then we may observe the
following from the definition itself.
$$
M\left(\alpha,{{p+1}\over2};A\right)=|A|^{-\alpha}\Gamma_{p}(\alpha),~{\rm
Re}(\alpha)>{{p-1}\over2}.\eqno(2.2)
$$
\noindent
But for $p=1$, that is , in the scalar variable case
the function $M(\cdot,\cdot;\cdot)$ is associated with a
Whittaker function. Hence we define a Whittaker function
in terms of $M(\cdot,\cdot;\cdot)$ as follows:

$$M(\mu,\nu;A)=|A|^{-{{\mu
+\nu}\over2}}\Gamma_{p}(\mu){\rm e}^{{1\over2}{\rm
tr}(A)}W_{{1\over2}(\nu -\mu),{1\over2}(\nu +\mu
-(p+1)/2)}(A)\eqno(2.3)
$$for $A=A'>0, {\rm Re}(\mu)>{{p-1}\over2}$ and no
restriction on $\nu$. In terms of the integral
representation we have
$$
\eqalignno{\int_{Z>0}|Z|^{\mu -{{p+1}\over2}}&|I+Z|^{\nu
-{{p+1}\over2}}{\rm e}^{-{\rm tr}(AZ)}{\rm d}Z\cr
&=|A|^{-{{\mu +\nu}\over2}}\Gamma_{p}(\mu){\rm
e}^{{1\over2}{\rm tr}(A)}W_{{1\over2}(\nu
-\mu),{1\over2}(\nu +\mu -(p+1)/2)}(A).&(2.4)\cr}
$$By using (2.4) one can establish a number of  results
on Whittaker function which will generalize the
corresponding univariate results to the matrix-variate
case. For the sake of illustration one such result will
 be listed here, before going back to the density in (1.10).
  This result follows from the above definition itself.
\vskip.6cm
\noindent
{\bf Theorem 2.1.}{\hskip.5cm}{\it For } ${\rm Re}(\beta
-\alpha)>{{p-3}\over4}, A=A'>0$,
$$
\eqalignno{\int_{Z>0}&|Z|^{\beta -\alpha
-{{p+1}\over4}}|I+Z|^{\alpha
+\beta -{{p+1}\over4}}{\rm
e}^{-{\rm tr}(AZ)}{\rm d}Z\cr
&=|A|^{-\beta -{{p+1}\over4}}\Gamma_{p}\left(\beta
-\alpha +{{p+1}\over4}\right){\rm e}^{{1\over2}{\rm
tr}(A)}W_{\alpha,\beta}(A).\cr}
$$This follows from the definition itself by taking $\alpha
={1\over2}(\nu -\mu)\hbox{ and }\beta ={1\over2}\left(\nu
+\mu -{{p+1}\over2}\right)$.
\vskip.2cm
Compare
the integral in (1.10) with the result in Theorem 2.1 to obtain
$$\eqalignno{g_1(U)&=\delta
|B_1+B_2|^{-{{\alpha_1+\alpha_2}\over2}}\Gamma_{p}(\alpha_2)h_1(U)\cr
\noalign{\hbox{where}}
h_1(U)&={\rm e}^{{1\over2}{\rm
tr}((B_2-B_1)U)}|U|^{{{\alpha_1+\alpha_2}\over2}-{{p+1}\over2}}W_{{1\over2}(
\alpha_1-\alpha_2),{1\over2}(\alpha_1+\alpha_2-(p+1)/2)}\left((B_1+B_2)^{1\o
ver2}U(B_1+B_2)^{1\over2}\right).&(2.5)\cr}
$$Now for $U<0, V>0$ let
$$g_2(U)=\int_{V>0}g(U,V){\rm d}V.$$Write $W=-U$ so that
$W=W'>0$. After taking out $W$  make the change
$Z=W^{-{1\over2}}VW^{-{1\over2}}$ for fixed $W$ to obtain
$$
\eqalignno{g_2(U)&=\delta |W|^{\alpha_1+\alpha_2-{{p+1}\over2}}{\rm
e}^{-{\rm tr}(B_2W)}\cr
&\times\>\int_{Z>0}|Z|^{\alpha_1-{{p+1}\over2}}|I+Z|^{\alpha_2-{{p+1}\over2}
}{\rm
e}^{-{\rm
tr}\left(W^{1\over2}(B_1+B_2)W^{1\over2}Z\right)}{\rm
d}Z,~W=-U.\cr}
$$Now in the light of Theorem 2.1, we have, for
$W=-U$.
$$\eqalignno{g_2(U)&=\delta
|B_1+B_2|^{-{{\alpha_1+\alpha_2}\over2}}\Gamma_{p}(\alpha_1)h_2(W)\cr
\noalign{\hbox{where}}
h_2(W)&=|W|^{{{\alpha_1+\alpha_2}\over2}-{{p+1}\over2}}{\rm
e}^{{1\over2}{\rm
tr}((B_1-B_2)W)}W_{{1\over2}(\alpha_2-\alpha_1),{1\over2}(\alpha_1+\alpha_2-
(p+1)/2)}\left((B_1+B_2)^{1\over2}W(B_1+B_2)^{1\over2}\right).&(2.6)\cr}
$$Hence the density of $Y=X_1-X_2$, in the oriented case discussed in (1.10),
 is given by the following
\vskip.5cm
\noindent
{\bf Theorem 2.2.}
$$g(Y)=\cases{c^{-1}h_1(Y),~Y>0\cr
c^{-1}h_2(W),~W=-Y>0\cr
0,\hbox{  elsewhere,}\cr}\eqno(2.7)
$$where $h_1(Y)$ and $h_2(W)$ are given in (2.5) and
(2.6) respectively, $c=c_1+c_2$, with
$$c_1=\int_{Y>0}h_1(Y){\rm d}Y\hbox{  and
}c_2=\int_{W>0}h_2(W){\rm d}W.
$$The density in (2.7) generalizes the univariate case
Mathai-1,4 of  Mathai (1993a, p.34).
\vskip.2cm
Integrals of the type in (1.10) appear in a large variety of statistical
 distribution problems connectd with sum, difference and linear functions of
 positve random variables, especially of the gamma type, or when dealing with
 the upper part of the incomplete gamma function. The corresponding
 matrix-variate analogues can be handled by using the definition of Whittaker
 function given above. When dealing with the distribution of a linear function
 of three independent real gamma variables one has to evaluate a double integral
 where the kernel part is of the type in (1.10). This leads to integrals
 involving Whittaker functions, see Mathai and Pederzoli (1996).
(Distributions of linear functions of gamma
 variables are discussed in Mathai and Provost (1992)). Several results on
 integrals involving Whittaker functions of matrix argument are given by the
 author and his co-workers recently. For the sake of illustration one such
 result will be given here as a theorem.
\vskip.5cm
\noindent
Theorem 2.3.\hskip.5cm {\it For $A=A'>0,\ B=B'>0,\  {\rm Re}\ (\gamma \pm
\beta ) > {p-3\over 4},\  {\rm Re}\ (\beta -\alpha )>{p-3\over 4},\ Z=Z'>0$,
$$\eqalignno{&\int_{Z>0} \vert Z\vert^{\gamma-{p+1\over 2}} {\rm e}^{-{\rm
tr}  (AZ)}  W_{\alpha,\beta} (BZ) {\rm d} Z =\left\vert A+{1\over 2}
B\right\vert^{-\left( \gamma +\beta +{p+1\over 4}\right)}\cr &\times
{\Gamma_p\left(\gamma +\beta +{p+1\over 4}\right) \Gamma_p \left(\gamma
-\beta +{p+1\over 4}\right)\over \Gamma_p \left(\gamma -\alpha +{p+1\over 2}
\right)}\cr &\times{}_2F_1 \bigg( {p+1\over 4} +\beta -\alpha,\
 {p+1\over 4}+\gamma +\beta ;\cr
&{p+1\over 2}+\gamma -\alpha ; I-B^{1\over 2}
\left(A+{1\over 2} B \right)^{-1} B^{1\over 2} \bigg) &(2.8)\cr
&\hskip 1cm\hbox {for}\  0<B^{1\over 2} \left( A+{1\over 2} B\right)^{-1}
B^{1\over 2} < I\ \hbox{or}\  2B^{-{1\over 2}} A B^{-{1\over 2}} >I}$$
$$\eqalignno{&=\vert B\vert^{-\left(\gamma +\beta +{p+1\over 4} \right)}
{\Gamma_p  \left({p+1\over 4} +\gamma +\beta \right) \Gamma_p \left({p+1\over
4} +\gamma -\beta\right) \over \Gamma_p \left({p+1\over 2} +\gamma -\alpha
\right)}\cr &\times{}_2F_1 \left({p+1\over 4} +\gamma -\beta , {p+1\over 4}
+\gamma +\beta ;  {p+1\over 2} +\gamma -\alpha; - C,\right )&(2.9)\cr
&\hskip7cm \hbox { for}\ \Vert C\Vert <1}$$ where
$$C=B^{-{1\over 2}} A B^{-{1\over 2}}-{1\over 2} I$$
and $\Vert (\cdot )\Vert$ denotes a norm of $(\cdot )$. A sufficient condition
is that
$$0<B^{-{1\over 2}}AB^{-{1\over 2}}-{1\over 2} I<I\ \hbox {or}\  0<{1\over 2}
I-B^{-{1\over 2}} AB^{-{1\over 2}} <I.$$}
\vskip.5cm
\noindent
 Proof.\hskip.5cm Replace
$W_{\alpha,\beta}(BX)$ by its integral representation to obtain
$$\eqalignno{\hbox{left side} &= \int_{X>0}\vert X\vert^{\gamma-{p+1\over 2}}
{\rm e}^{-{\rm tr} (AX)} W_{\alpha,\beta}(BX) {\rm d}X\cr
&=\int_{X>0} {\vert X\vert^{\gamma +\beta +{p+1\over 4}-{p+1\over 2}}\over
\Gamma_p \left({p+1\over 4}+\beta -\alpha \right)} {\rm e}^{-{\rm tr}
(AX)-{1\over 2} {\rm tr} (BX)} \cr
& \times \int_{Z>0} \vert Z\vert^{\beta -\alpha -{p+1\over 4}}
\vert I+Z\vert^{\alpha+\beta-{p+1\over 4}}
{\rm e}^{-{\rm tr} (BX Z)} {\rm d}Z {\rm d}X.}$$
The $X$-integral can be evaluated  by using a gamma integral to get
$$\eqalignno{\int_{X>0} &\vert X\vert^{\gamma +\beta +{p+1\over 4}-{p+1\over
2}}  {\rm e}^{-{\rm tr}  \left[\left((A+{1\over 2} B)+BZ\right)X\right]} {\rm
d}X \cr &= \Gamma_p \left(\gamma +\beta
+{p+1\over 4} \right) \left\vert \left(A+{1\over 2} B\right ) +BZ
\right\vert^{-\left(\gamma +\beta +{p+1\over 4}\right)}\cr
&\hskip2cm \hbox{for}\ {\rm Re} (\gamma +\beta )>{p-3\over 4}.&(2.10)}$$
The $Z$-integral, denoted by $f$, is given by
$$\eqalignno{f &= \int_{Z>0} \vert Z\vert^{\beta -\alpha -{p+1\over 4}}
\vert I+Z\vert^{\alpha +\beta -{p+1\over 4}} \left\vert \left (A+{1\over
2}B\right) +BZ\right\vert^{-\left( \gamma +\beta +{p+1\over 4} \right)} {\rm
d} Z\cr &= \left\vert A+{1\over 2} B\right\vert^{-\left(\gamma +\beta
+{p+1\over 4}\right)} \int_{Z>0} \vert Z\vert^{\beta-\alpha -{p+1\over 4}}
\vert I +Z\vert^{\alpha +\beta-{p+1\over 4}}\cr
&\times \left\vert I+B^{1\over 2} \left( A+{1\over 2} B\right)^{-1}
B^{1\over 2} Z \right\vert^{-\left(\gamma +\beta +{p+1\over 4}\right)}{\rm d}
Z.}$$ Note that in the determinant we can also replace
$\left( A+{1\over 2} B\right)^{-1} B $ by $B^{1\over 2} \left( A+{1\over 2}
B\right)^{-1} B^{1\over 2}$. Write $Z= (I-U)^{-{1\over 2}}U(I-U)^{-{1\over 2}}
\Rightarrow I+Z^{-1} = U^{-1}  \Rightarrow \vert
Z\vert ^{-(p+1)} {\rm d}Z = \vert U\vert ^{-(p+1)} {\rm d}U \Rightarrow {\rm
d}Z =\vert I-U\vert ^{-(p+1)} {\rm d}U$ and $0<U<I$. Then
$$\eqalignno{f=&\left\vert A +{1\over 2} B\right\vert^{-(\gamma +\beta
+{p+1\over 4}} \int_{0<U<I} \vert U\vert^{\beta -\alpha -{p+1\over 4}} \vert
I-U\vert^{\gamma -\beta -{p+1\over 4}}\cr  &\times \left\vert
I-\left(I-B^{1\over 2}\left(A+{1\over 2} B\right)^{-1} B^{1\over
2}\right)U\right\vert^{-\left(\gamma +\beta +{p+1\over 4}\right)} {\rm d}U.}$$
Now we can write the integral as a\   ${}_2F_1$ by using Mathai (1993a, p.179).
  Then $$\eqalignno{f&=\vert A+{1\over 2} B\vert^{-\left(\gamma +\beta
+{p+1\over 4}\right)}\cr &\times{\Gamma_p\left(\beta -\alpha +{p+1\over 4}
\right) \Gamma_p \left(\gamma - \beta +{p+1\over 4}\right)\over
\Gamma_p\left(\gamma -\alpha +{p+1\over 2}\right)}\cr
&\times{}_2F_1\bigg({p+1\over 4}+\beta -\alpha, {p+1\over 4} +\gamma
+\beta;\cr &{p+1\over 2}+\gamma -\alpha ; I-B^{1\over 2} \left( A+ {1\over 2}
B\right)^{-1} B^{1\over 2} \bigg)}$$ for $0<I-B^{1\over 2} \left( A+ {1\over
2} B\right)^{-1} B^{1\over 2} < I  \Rightarrow 0<B^{1\over 2} \left(
A+{1\over 2} B\right)^{-1} B^{1\over 2} < I \Rightarrow 2B^{-{1\over 2}}
AB^{-{1\over 2}} >I$. Substituting back, one gamma is cancelled, we get
(2.8). But observe that
$$\eqalignno{\left\vert \left( A+{1\over 2} B\right) + BZ\right\vert
&=\left\vert \left( A-{1\over 2} B\right) +B\left( I-Z\right)\right\vert\cr
 &= \vert B\vert\  \vert I-Z\vert\  \left\vert I +\left( B^{-{1\over 2}}
 AB^{-{1\over 2}}-{1\over 2} I \right) \left( I+Z\right) ^{-1} \right\vert.} $$
Now

$$\eqalignno{f& = \vert B\vert^{-\left( \gamma +\beta +{p+1\over 4}\right)}
\int_{Z>0}\vert Z\vert^{\beta -\alpha -{p+1\over 4}} \vert I+Z\vert^{\alpha
-\gamma -{p+1\over 2}} \cr &\times \left\vert I+\left( B^{-{1\over 2}}
AB^{-{1\over 2}}-{1\over 2} I\right)  \left( I + Z\right)^{-1} \right\vert
^{-\left(\gamma +\beta +{p+1\over 4}\right)} {\rm d} Z.}$$
Put $U=(I+Z)^{-1} \Rightarrow {\rm d} Z=\vert U\vert^{-(p+1)} {\rm d} U,\
\vert Z\vert = \vert U \vert^{-1} \vert I-U\vert, \vert I +Z\vert =\vert
U\vert ^{-1}, 0<U<I$. Then
$$f=\vert B\vert^{-\left(\gamma +\beta +{p+1\over 4}\right)} \int_{0<U<I}
\vert U\vert^{\gamma -\beta -{p+1\over 4}} \vert I - U\vert^{\beta-\alpha
-{p+1\over 4}} \vert I+CU\vert^{-\left(\gamma +\beta +{p+1\over 4}\right)}$$
where
$$C=B^{-{1\over 2}} AB^{-{1\over 2}} -{1\over 2} I.$$
Evaluating this integral by using Mathai (1993a, p.179) one has
$$\eqalign{f&=\vert B\vert^{-\left(\gamma +\beta +{p+1\over 4}\right)}
{\Gamma_p\left({p+1\over 4}+\gamma -\beta \right) \Gamma_p \left( {p+1\over
4}+\beta -\alpha \right) \over \Gamma_p \left({p+1\over 2} +\gamma -\alpha
\right)}\cr &\times{}_2F_1 \left({p+1\over 4}+\gamma -\beta, \gamma +\beta
+{p+1\over 4};  {p+1\over 2} +\gamma -\alpha ; - C\right),\cr
 &\hbox {for}\ \Vert C\Vert <1,\ {\rm Re}\ (\gamma -\beta )>{p-3\over 4},
{\rm Re}\ (\beta
-\alpha ) > {p-3\over 4}.}$$ \noindent
Substituting back we get (2.9), noting that one gamma is cancelled.
\vskip.5cm
\noindent
3.{\hskip.5cm}RESIDUAL VARIABLES IN THE COMPLEX CASE
\vskip.5cm
Let $\tilde{X}_1$ and $\tilde{X}_2$ be $p\times p$ hermitian positive definite
 matrix random variables defined in the complex field. Such a matrix
 $\tilde{X}_j$ can be expressed in the following form:
$$\tilde{X}_j=X_{j1}+iX_{j2}
$$where $i=\sqrt{-1}$, $X_{j1}$ and $X_{j2}$ are matrices with real elements
 such that $X_{j1}=X_{j2}'>0$ and $X_{j2}'=-X_{j2}$. That is, $X_{j1}$ is
 symmetric positive definite and $X_{j2}$ is skew symmetric. All the matrices
 appearing in this section are $p\times p$ hermitian positive definite unless
 stated otherwise.
\vskip.2cm
If the growth and decay matrices in the complex field are respectively
 $\tilde{X}_1$ and $\tilde{X}_2$ then the residual variable is
$$\tilde{Y}=\tilde{X}_1-\tilde{X}_2.\eqno(3.1)
$$Let $\tilde{X}_1$ and $\tilde{X}_2$ be independently distributed
 matrix-variate complex gamma variables with the densities
$$\tilde{g}_j(\tilde{X}_j)={{|{\rm
det}(\tilde{B}_j)|^{\alpha_j}}\over{\tilde{\Gamma}_p(\alpha_j)}}|{\rm
det}(\tilde{X}_j)|^{\alpha_j-p}{\rm e}^{-{\rm
tr}(\tilde{B}_j\tilde{X}_j)},~\tilde{X}_j=\tilde{X}_j^{*}>0,~{\rm
Re}(\alpha_j)>p-1\eqno(3.2)
$$ and $\tilde{g}_j(\tilde{X}_j)=0$ elsewhere, where
 $(\cdot)^{*}$ denotes the conjugate transpose of $(\cdot)$,
 $\tilde{X}_j=\tilde{X}_j^{*}>0$ implies that $\tilde{X}_j$ is hermitian
 positive definite, ${\rm det}(\cdot)$ denotes the determinant of $(\cdot)$,
 $|{\rm det}(\cdot)|$ gives the absolute value of the determinant and
 $\tilde{\Gamma}_p(\cdot)$ is the matridx-variate gamma in the complex case,
 given by
$$\eqalignno{\tilde{\Gamma}_p(\alpha)&=\int_{\tilde{X}=\tilde{X}^{*}>0}|{\rm
det}(\tilde{X})|^{\alpha -p}{\rm e}^{-{\rm tr}(\tilde{X})}{\rm
d}\tilde{X}&(3.3)\cr
&=\pi^{p(p-1)/2}\Gamma(\alpha)\Gamma(\alpha -1)...\Gamma(\alpha -p+1),~{\rm
Re}(\alpha )>p-1.\cr}
$$Quadratic forms in complex Gaussian random vectors, giving rise to Whittaker
 variables, are used in a wide variety of problems in communication and
 engineering areas, see for example, Hirasawa (1988), Divsalar et al.(1990),
 and Biyari and Lindsey (1991).
\vskip.2cm
Our aim here is to look at the residual variable in (3.1) when the input and
 output variables are independently distributed with the densities in (3.2).
 Then proceeding as in (1.7) to (1.10) we note that in order to compute the
 density of the residual variable one needs an integral to be evaluated,
 corresponding to the one in (1.10). This requires the definition of a Whittaker
 function of matrix argument in the complex case.
 Whittaker function of matrix argument in the
 complex case, denoted by $\tilde{W}_{\cdot,\cdot}(\cdot)$, will be defined
as that
 symmetric function, symmetric in the sense
$$\tilde{W}_{\alpha,\beta}(\tilde{A}\tilde{B})=\tilde{W}_{\alpha,\beta}
(\tilde{B}\tilde{A})
$$for $p\times p$ hermitian positive definite matrices, having the following
 integral representaton.
$$\int_{\tilde{Z}=\tilde{Z}^{*}>0}|{\rm det}(\tilde{Z})|^{\beta -\alpha
-{p\over2}}|{\rm det}(I+\tilde{Z})|^{\beta +\alpha -{p\over2}}~{\rm
e}^{-{\rm tr}(\tilde{A}\tilde{Z})}{\rm d}\tilde{Z}$$
$$=|{\rm
det}(\tilde{A})|^{-\beta-{p\over2}}~\tilde{\Gamma}_{p}\left(\beta-\alpha
+{p\over2}\right){\rm e}^{{1\over2}{\rm
tr}(\tilde{A})}\tilde{W}_{\alpha,\beta}(\tilde{A})\eqno(3.4)
$$for ${\rm Re}(\beta -\alpha)>{p\over2}-1$.
\vskip.2cm
As a consequence of the definition itself we can have the following result
 which will be stated as a theorem.
\vskip.7cm
\noindent
{\bf Theorem 3.1.}{\hskip.5cm}{\it For $\tilde{A}=\tilde{A}^{*}>0$}
$$\eqalignno{\int_{\tilde{Z}=\tilde{Z}^{*}>0}&|{\rm
det}(I+\tilde{Z})|^{-\alpha}{\rm e}^{-{\rm tr}(\tilde{A}\tilde{Z})}{\rm
d}\tilde{Z}\cr
&=|{\rm det}(\tilde{A})|^{{\alpha\over2}-p}~\tilde{\Gamma}_{p}(p){\rm
e}^{{1\over2}{\rm
tr}(\tilde{A})}\tilde{W}_{-{\alpha\over2},{1\over2}(-\alpha
+p)}(\tilde{A}).\cr}
$$One can establish a number of results on Whittaker function in the complex
 case which will be useful in working out the distributions of sums, differences
 and linear functions of matrix gamma  variables in the complex case, evaluating
 probabilities associated with these distributions, evaluating the incomplete
 gamma, incomplete type-2 beta and other related integrals, and related
problems.
 In order to illustrate the techniques one more result will be given here.
\vskip.5cm
\noindent
{\bf Theorem 3.2.}{\hskip.5cm}{\it For $\tilde{B}=\tilde{B}^{*}>0,~
\tilde{U}=\tilde{U}^{*}>0,~ \tilde{M}=\tilde{M}^{*}>0,~ {\rm
Re}(q)>{p\over2}$},
$$\eqalignno{\int_{\tilde{X}>\tilde{U}}&|{\rm
det}(\tilde{X}+\tilde{B})|^{2\alpha -p}|{\rm
det}(\tilde{X}-\tilde{U})|^{2q-p}{\rm e}^{-{\rm
tr}(\tilde{M}\tilde{X})}{\rm d}\tilde{X}\cr
&=\tilde{\Gamma}_{p}(2q)|{\rm det}(\tilde{U}+\tilde{B})|^{\alpha +q-p}|{\rm
det}(\tilde{M})|^{-(\alpha +q)}\cr
&\times {\rm e}^{{1\over2}{\rm
tr}[(\tilde{B}-\tilde{U})\tilde{M}]}\tilde{W}_{(\alpha -q),(\alpha
+q-{p\over2})}(\tilde{T}),\cr
\tilde{T}&=(\tilde{U}+\tilde{B})^{{1\over2}}\tilde{M}(\tilde{U}+\tilde{B})^{
{1\over2}}.\cr}$$
\vskip.3cm
\noindent
Proof.{\hskip.5cm}Consider the following transformations.
$$\eqalignno{\tilde{Y}&=\tilde{X}-\tilde{U},~\tilde{Z}=(\tilde{B}+\tilde{U})
^{-{1\over2}}\tilde{Y}(\tilde{B}+\tilde{U})^{-{1\over2}}\cr
&\Rightarrow{\rm d}\tilde{Z}=|{\rm det}(\tilde{B}+\tilde{U})|^{-p}{\rm
d}\tilde{Y}.\cr}
$$Now the integral on the left side, denoted by $f$, is given by the following:
$$\eqalignno{f&={\rm e}^{-{\rm tr}(\tilde{M}\tilde{U})}|{\rm
det}(\tilde{B}+\tilde{U})|^{2\alpha +2q-p}\cr
&\times\int_{\tilde{Z}=\tilde{Z}^{*}>0}|{\rm det}(\tilde{Z})|^{2q-p}|{\rm
det}(I+\tilde{Z})|^{2\alpha -p}\cr
&\times\exp\left\{-{\rm
tr}\left[\tilde{M}(\tilde{B}+\tilde{U})^{1\over2}\tilde{Z}(\tilde{B}+\tilde{
U})^{1\over2}\right]\right\}{\rm d}\tilde{Z}\cr}
$$for ${\rm Re}(2q)>p-1$. Now interpreting with the help of the definition of a
 Whittaker function
 the result follows.
\vskip.5cm
\noindent
\centerline{REFERENCES}
\vskip.5cm
\noindent
Biyari, K.H. and Lindsey, W.C. (1991). Statistical distributions of hermitian
 quadratic forms in complex Gaussian variables. {\it IEEE Transactions on
 Information Theory}, {\bf 39(3)}, 1076--1082.
\vskip.3cm
\noindent
Divsalar, E. et al. (1990). The performance of trellis-coded MDPSK with multiple
 symbol detection. {\it IEEE Transactions on Communcations},  {\bf 38(9)},
1391--1403.
\vskip.3cm
\noindent
Haubold, H.J. and Mathai, A.M. (1995). A heuristic remark on the periodic
 variation in the number of solar neutrinos detected on earth. {\it Astrophysics
 and Space Science}, {\bf 228}, 113--134.
\vskip.3cm
\noindent
Hirasawa, K. (1988). The application of a biquadratic programming method of
 phase-only optimization of antenna arrays. {\it IEEE Transactions on Antenna
 Propagation}, {\bf 36(11)}, 1545--1550.
\vskip.3cm
\noindent
Mathai, A.M. (1993a). {\it A Handbook of Generalized Special Functions for
 Statistical and Physical Sciences}, Oxford University Press, Oxford.
\vskip.3cm
\noindent
Mathai, A.M. (1993b). The residual effect of a growth-decay mechanism and the
 distributions of covariance structures. {\it The Canadian Journal of
 Statistics}, {\bf 21(3)}, 277--283.
\vskip.3cm
\noindent
Mathai, A.M. and Pederzoli, G. (1996). Some transformations for functions of
 matrix arguments. {\it Indian Journal of Pure and Applied Mathematics},
 {\bf 27(3)}, 277-284.
\vskip.3cm
\noindent
Mathai, A.M.  and Provost, S.B. (1992). {\it Quadratic Forms in Random
Variables:
 Theory and Applications}, Marcel Dekker, New York.
\vskip.3cm
\noindent
Mathai, A.M., Provost, S.B. and Hayakawa, T. (1995). {\it Bilinear Forms  and
 Zonal Polynomials}, Springer-Verlag, (Lecure Notes in Statistics), New York.
\vskip.3cm
\noindent
Rao, C.R. (1973). {\it Linear Statistical Inference and Its Applications},
 Wiley, New York.
\bye